\documentclass[11pt]{amsart}
\usepackage[margin=1in]{geometry}
\usepackage{amsmath,amssymb,graphicx,booktabs,hyperref}
\newtheorem{theorem}{Theorem}
\newtheorem{proposition}[theorem]{Proposition}
\newtheorem{lemma}[theorem]{Lemma}
\newtheorem{corollary}[theorem]{Corollary}
\newtheorem{conjecture}{Conjecture}
\theoremstyle{remark}
\newtheorem{remark}[theorem]{Remark}
\newcommand{\LT}{\widetilde L}
\newcommand{\cL}{\mathcal L}
\newcommand{\C}{\mathbb C}
\newcommand{\R}{\mathbb R}
\def\x#1{}
\def\xyz{\hspace*{.07em}}

\def\xz{\hspace*{-.07em}}

\DeclareMathOperator{\Ree}{Re}
\DeclareMathOperator{\Imm}{Im}

\title[How $\Lambda_n$ leaves the polygon]{How the eigenvalue region of standardized Laplacians\\ leaves the 2005 polygon}
\author{Pavel Chebotarev}
\address{Department of Mathematics, Technion---Israel Institute of Technology, Haifa 3200003, Israel}
\email{pavel4e@technion.ac.il}
\date{DRAFT, 31 October 2026}

\begin{document}
\begin{abstract}
Let $\Lambda_n$ be the set of all eigenvalues of $n\times n$ standardized Laplacian
matrices, and let $S_n\subseteq\Lambda_n$ be the polygon introduced in 2005, whose upper
boundary is the chain $0=\lambda_0,\lambda_1,\dots,$ $\lambda_{n-1}=1$. It is now known
that $\Lambda_n=S_n,$ and that $\Lambda_n$ is convex, if and only if $n\le4$. We study how
$\Lambda_n$ leaves $S_n$ near an edge $[\lambda_m,\lambda_{m+1}]$ of $S_n$ with
$2(m+1)<n$.
First, we prove a switching law for matrices having a boundary point of $\Lambda_n$ as
a simple eigenvalue: an explicit switching function on the arcs, built from the left and
right eigenvectors and the outward normal direction, determines the weights. Every arc
on which it is positive has the maximal weight, every arc on which it is negative has
weight zero, and intermediate weights can occur only on arcs where it vanishes.
Second, the point $\lambda_m+\tau(\lambda_{m+1}-\lambda_m)$ of the edge is realized by
the circulant digraph with arcs $k\to k+s,$ $1\le s\le m,$ and arcs $k\to k+m+1$ of
weight $\tau$. We perturb the weights of the latter arcs. At second order, the first
Fourier harmonic of the perturbation moves the eigenvalue beyond the line of the edge if
and only if
\[
\tau<\tau^*_m=\frac{\sin(\pi m/n)}{2\cos(\pi/n)\,\sin(\pi(m+1)/n)}.
\]
For even $n$ and $m,$ the point of the edge with $\tau=\tau^*_m$ is the vertex of the
inscribed polygon $S_{n/2}$ that lies on this edge. An explicit one-parameter family of
digraphs, in which the arcs $k\to k+m+1$ are present only for even $k,$ has an
eigenvalue path that leaves $S_n$ at $\lambda_m$ and, in all computed cases, ends
exactly at this vertex.
Third, we give a conjectural, numerically supported description of $\partial\Lambda_5$.
Beyond each of the four edges of $S_5$ not incident to $0$ or $1,$ $\Lambda_5$ forms a
cap bounded by four algebraic arcs; it leaves the edge tangentially at the endpoint
adjacent to $0$ or $1$ and rejoins it at $\tau^*_1=(3-\sqrt5)/2,$ and the rest of the
edge lies on $\partial\Lambda_5$. Computations for $n\le8$ show that a cap can return to
its edge in three ways: smoothly at $\tau^*_m,$ at a corner at $\tau^*_m,$ or at a corner
beyond $\tau^*_m$.
\end{abstract}


\maketitle

\section{Introduction}

A \emph{standardized Laplacian matrix} of order $n$ is a real matrix $\LT=(\ell_{ij})$ with zero row sums and $-\frac1n\le\ell_{ij}\le0$ for $j\ne i$. Let $\cL_n$ be the set of these matrices and $\Lambda_n=\bigcup_{\LT\in\cL_n}\sigma(\LT)$. Put $\mu=e^{-2\pi i/n},$ $g_s=\frac1n(1-\mu^s),$ and $\lambda_k=g_1+\dots+g_k,$ $0\le k\le n-1$. The polygon $S_n=\operatorname{conv}\{\lambda_k,\bar\lambda_k\}$ is contained in $\Lambda_n$ \cite[Theorem~6]{AgaChe05LAA}, and in 2005 it was conjectured that $\Lambda_n=S_n$ \cite{AgaChe05LAA}. The conjecture holds for $n\le3$ \cite{AgaChe05LAA} and for $n=4$ \cite{Liu026,Marshall026,Liu026b}, and fails for every $n\ge5$ \cite{Che26arXi}. Moreover, every edge of $S_n$ not incident to $0$ or $1$ contains interior points of $\Lambda_n,$ except possibly the two horizontal edges for even $n$ \cite[Corollary~16, Remark~9]{Che26arXi}, and $\Lambda_n$ is not convex for $n\ge5$ \cite{Liu026}, \cite[Theorem~4]{Liu026b}. The set of points of $\Lambda_n$ lying beyond the line of an edge of $S_n$ whose orthogonal projection onto that line falls in the edge we call the \emph{cap\/} over that edge.

The results for $n\ge5$ are mostly local: they come from second-order perturbations of circulant matrices at the vertices of $S_n$. Here we ask what the boundary of $\Lambda_n$ actually looks like near an edge. Our answers are of three kinds.

\emph{Theorems valid for all $n$.} A switching law (Theorem~\ref{t:switch}) gives necessary conditions on matrices realizing boundary points.
Wherever the switching function is nonzero, the arc weights are extremal ($0$ or maximal; a bang--bang structure). 
Intermediate weights, i.e., weights strictly between $0$ and $\frac1n,$ can occur only where it vanishes. They satisfy explicit fold equations. 
For a fixed combinatorial type, these equations generically constrain the candidate boundary to an algebraic curve.
A reversal symmetry (Lemma~\ref{l:rev}) underlies the observed symmetry, under a reflection of the cycle, of the chord weights realizing the caps ({\em chords\/} are the arcs $k\to k+m+1$ of the digraph). Along an edge, the second-order coefficient of the first harmonic of chord perturbations changes sign at the explicit threshold $\tau^*_m$ (Theorem~\ref{t:thresh}). For even $n$ and $m,$ alternating chords give an explicit arc that ends at a vertex of $S_{n/2}$ (Proposition~\ref{p:alt}).

\emph{A conjectural complete description of the boundary for $n=5$} (Conjecture~\ref{c:five}). Each piece is given by explicit equations, the transition points are computed to many digits, and the description is supported by a global numerical search.

\emph{Observations for $n\le8$.} We describe three different ways in which a cap returns to the edge. We also describe the growth of palindromic chord blocks, whose changes of parity are accompanied by small corners inside the caps.

Throughout, the vertices of the digraphs are labeled $0,\dots,n-1,$ indices are taken modulo $n,$ $Q$ is the cyclic permutation matrix with $(Qx)_k=x_{k+1},$ $f_k=(\mu^{kl})_{l=0}^{n-1},$ and $D_\omega=\operatorname{diag}(\omega_0,\dots,\omega_{n-1})$. By \cite[Corollary~22]{Che26arXi}, $\Lambda_n$ is star-shaped with respect to every point of $[0,1]$.

\section{The switching law}

Write $a_{ij}=-\ell_{ij}\in\bigl[0,\frac1n\bigr]$. If $\lambda$ is a simple eigenvalue of $\LT$ with right and left eigenvectors $v$ and $w,$ then
\begin{equation}\label{e:grad}
\frac{\partial\lambda}{\partial a_{ij}}=\frac{w_i(v_i-v_j)}{w^{\mathsf T}v},
\end{equation}
since changing $a_{ij}$ by $\delta$ adds $\delta$ to the entry $(i,i)$ and $-\delta$ to the entry $(i,j)$.

Let $\nu$ be a unit complex number. We say that $\nu$ is an \emph{outward normal} of $\Lambda_n$ at $\lambda\in\Lambda_n$ if for every closed cone $K\subset\{u:\Ree(\bar\nu u)>0\}\cup\{0\}$ there is $\varepsilon>0$ such that $(\lambda+K)\cap\Lambda_n\cap\{|z-\lambda|<\varepsilon\}=\{\lambda\}$. This holds at every point of a $C^1$ arc of $\partial\Lambda_n$ near which $\Lambda_n$ lies on one side of the arc (with its outward normal) and at every point where $\Ree(\bar\nu z)$ attains a local maximum on $\Lambda_n$.

\begin{theorem}[switching law]\label{t:switch}
Let $\nu$ be an outward normal of $\Lambda_n$ at $\lambda,$ and let $\lambda$ be a simple eigenvalue of $\LT\in\cL_n$ with right and left eigenvectors $v,w$. Put
\[
G_{ij}=\Ree\Bigl(\bar\nu\,\frac{w_i(v_i-v_j)}{w^{\mathsf T}v}\Bigr),\qquad i\ne j.
\]
Then $a_{ij}=\frac1n$ whenever $G_{ij}>0,$ and $a_{ij}=0$ whenever $G_{ij}<0$. In particular\textup, $a_{ij}\in\bigl(0,\frac1n\bigr)$ only if $G_{ij}=0$.
\end{theorem}

\begin{proof}
Suppose $G_{ij}>0$ and $a_{ij}<\frac1n$. Increasing $a_{ij}$ by $t\ge0$ keeps the matrix in $\cL_n$ for small $t$. By \eqref{e:grad}, the eigenvalue becomes $\lambda(t)=\lambda+t\delta+O(t^2)$ with $\Ree(\bar\nu\delta)=G_{ij}>0$. For a closed cone $K$ around $\delta$ inside $\{\Ree(\bar\nu u)>0\},$ we have $\lambda(t)\in\lambda+K$ for small $t>0,$ which contradicts the definition of an outward normal. The case $G_{ij}<0,$ $a_{ij}>0$ is symmetric.
\end{proof}

\begin{remark}
The simplicity assumption in Theorem~\ref{t:switch} is needed for the eigenvalue derivative formula \eqref{e:grad}.
At all points listed in Table~\ref{tab:five}, the relevant eigenvalues are simple and well conditioned.
Whether a boundary point of $\Lambda_n$ can be realized only as a multiple eigenvalue remains open.
\end{remark}

\begin{corollary}[fold equations]\label{c:fold}
Under the hypotheses of Theorem~\textup{\ref{t:switch},} let $\mathcal I$ be the set of arcs with intermediate weights\textup, and let $P(z;a)=\det(zI-\LT(a))$. Then $\Ree(\bar\nu\,\partial\lambda/\partial a_e)=0$ for all $e\in\mathcal I$. Equivalently\textup,
\[
\Imm\bigl(P_{a_e}\,\overline{P_{a_{e'}}}\bigr)=0\quad(e,e'\in\mathcal I),\qquad \bar\nu\,P_{a_e}/P_z\in i\xyz\R .
\]
\end{corollary}

Thus, along an arc of $\partial\Lambda_n$ on which the realizing matrices keep their \emph{type} (the sets of zero, maximal and intermediate weights), the boundary is cut out by $P=0$ and $|\mathcal I|-1$ real fold equations. These are $|\mathcal I|+1$ real equations in $|\mathcal I|+2$ real unknowns, so generically they define a curve, whose projection to the $\lambda$-plane lies on an algebraic curve. For $|\mathcal I|=1$ there is no fold equation, and the arc is the eigenvalue curve of a one-parameter family, as for the Karpelevi\v{c} arcs \cite{Karpelevich51E}.

\begin{remark}[circulants]\label{r:circ}
For a circulant $\LT$ and $\lambda$ its eigenvalue at $f_1,$ we have $v=f_1,$ $w=\bar f_1,$ and $G_{ij}=\Ree(\bar\nu g_{j-i})$. These are the terms of the support function of the zonotope $S_n=\sum_s[0,g_s]$ \cite[Proposition~1]{Che26arXi}. For $\nu=\nu_m,$ the outward normal of $[\lambda_m,\lambda_{m+1}],$ the law makes all arcs of length $\le m$ maximal and all arcs of length $>m+1$ zero, while for arcs of length $m+1$ it gives $G=0$ and decides nothing. This degeneracy is the reason why the first-order analysis \cite[Proposition~13]{Che26arXi} cannot decide the fate of the edges, and the second order has to.
\end{remark}

\begin{remark}
In the geometric approach of Dmitriev and Dynkin \cite{DmitrievDynkin46E} and Karpelevi\v{c} \cite{Karpelevich51E}, a boundary point of the eigenvalue region is studied through the extremal position of the eigenvalue relative to the polygon spanned by the components of an eigenvector; in \cite{Che26arXi} this takes the form $\lambda\in\bigcap_i Z_i(v)$. For fixed $v,$ maximizing $\Ree(\bar\nu\lambda)$ over the weights is a linear program \cite[Proposition~17]{Che26arXi}, \cite[Lemma~16]{Liu026b}. Taking $\bar\nu\,w_i/(w^{\mathsf T}v)$ as the dual multipliers of the row equations satisfies dual feasibility in $\lambda,$ since $\sum_i\bar\nu\,w_iv_i/(w^{\mathsf T}v)=\bar\nu,$ and the $G_{ij}$ are exactly the corresponding reduced costs. Theorem~\ref{t:switch} is thus the complementary slackness condition of this program, with the left eigenvector supplying the dual multipliers.
\end{remark}

\section{A reversal symmetry}

For $1\le m\le n-2$ and $\omega\in\C^n$ put
\begin{equation}\label{e:Nm}
N_m(\omega)=\sum_{s=1}^m(I-Q^s)+D_\omega(I-Q^{m+1}),
\end{equation}
so that $\frac1nN_m(\omega)$ is the standardized Laplacian of the $m$th power of the Hamiltonian cycle with chords $k\to k+m+1$ of weights $\omega_k\in[0,1]$ \cite[(12)]{Che26arXi}.

\begin{lemma}[reversal]\label{l:rev}
Let $\omega^R_k=\omega_{-k}$. Then $N_m(\omega)$ and $N_m(\omega^R)$ have the same characteristic polynomial.
\end{lemma}

\begin{proof}
Write $I-Q^s=(I-Q)E_s$ with $E_s=\sum_{r<s}Q^r,$ and put $B=\sum_{s\le m}E_s$ and $C=E_{m+1}$. These are polynomials in $Q,$ and $N_m(\omega)=(B+D_\omega C)(I-Q)$. Let $R$ be the permutation matrix with $Re_k=e_{-k},$ so that $RQ^{\mathsf T}R=Q$ and $RD_\omega R=D_{\omega^R}$. Then
\[
RN_m(\omega)^{\mathsf T}R=(I-Q)(B+C D_{\omega^R}).
\]
If $C$ is invertible, this equals $C\,(I-Q)(C^{-1}B+D_{\omega^R})\,,$ which is similar to $(I-Q)(C^{-1}B+D_{\omega^R})C=(I-Q)(B+D_{\omega^R}C),$ because $B$ and $C$ commute. Since $XY$ and $YX$ have the same characteristic polynomial, the latter has the characteristic polynomial of $(B+D_{\omega^R}C)(I-Q)=N_m(\omega^R)$. If $C$ is singular, replace $C$ by $C+sI$ in both $N_m(\omega)$ and $N_m(\omega^R)$. The identity of characteristic polynomials holds for all but finitely many $s,$ hence for all $s,$ in particular for $s=0$.
\end{proof}

Consequently, if a chord-weight vector $\omega$ realizes a boundary point, so does every reflection $k\mapsto\omega_{c-k}$. In all our computations, the nonconstant weight vectors realizing points of $\partial\Lambda_n$ beyond the edges are invariant under such a reflection, with the center at a vertex or at the midpoint between two vertices (Section~\ref{s:five}).
At a palindromic weight vector $\omega$, if the relevant eigenvalue is simple, Lemma~\ref{l:rev} and cyclic relabeling imply $\partial\lambda/\partial\omega_k=\partial\lambda/\partial\omega_{c-k}.$ Thus the fold equations for mirror-image intermediate weights
coincide.

\section{A second-order threshold along an edge} 

Fix $1\le m$ with $2(m+1)<n,$ put $p=m+1,$ let $\nu_m=-e^{-i\pi p/n}$ be the outward normal of $[\lambda_m,\lambda_{m+1}],$ and let $d_m(z)=\Ree(\bar\nu_m(z-\lambda_m))$ \cite[Lemma~11]{Che26arXi}. For $\tau\in[0,1]$ the circulant $\frac1nN_m(\tau\mathbf 1)$ has the eigenvalue $\lambda_m+\tau g_p\in[\lambda_m,\lambda_{m+1}]$ at $f_1$. By Remark~\ref{r:circ}, only the chords $k\to k+p$ are undecided at first order there. We perturb them: $\omega=\tau\mathbf 1+t\eta,$ $\eta\in\R^n$.

\begin{theorem}\label{t:thresh}
For $\tau\in[0,1],$ the eigenvalue $\lambda(t)$ of $\frac1nN_m(\tau\mathbf 1+t\eta)$ with $\lambda(0)=\lambda_m+\tau g_p$ is analytic near $t=0$ and
\[
d_m(\lambda(t))=t^2\sum_{j=1}^{n-1}c_j(\tau)\,|\hat\eta(j)|^2+O(t^3),\qquad
c_j(\tau)=\frac2n\sin\frac{\pi p}{n}\,\Imm\frac{1-\mu^{p(j+1)}}{F(j+1)-F(1)+\tau\bigl(\mu^{p(j+1)}-\mu^p\bigr)},
\]
where $F(k)=\sum_{s=0}^m\mu^{sk}$ and $\hat\eta(j)=\frac1n\sum_k\eta_k\mu^{jk}$. Moreover\textup, $c_{n-1}\equiv0$ and
\begin{equation}\label{e:c1}
c_1(\tau)=\frac{2}{n}\,\frac{\sin\frac{2\pi p}{n}\,\sin\frac{\pi p}{n}}{|b(\tau)|^2}\,\bigl(\tau^*_m-\tau\bigr),\qquad
\tau^*_m=\frac{\sin\frac{\pi m}{n}}{2\cos\frac{\pi}{n}\,\sin\frac{\pi(m+1)}{n}}
=\frac{\sin\frac{\pi m}{n}}{\sin\frac{\pi m}{n}+\sin\frac{\pi(m+2)}{n}},
\end{equation}
with $b(\tau)=\frac{\mu^p-\mu}{1-\mu^2}-\tau\mu^p\ne0$. If $n$ and $m$ are even\textup, then $c_{n/2}(\tau)=\frac12c_1(\tau)$ for all $\tau$.
\end{theorem}

\begin{proof}
The expansion and the formula for $c_j$ come from the Schur-complement argument of \cite[Lemma~8, Proposition~14]{Che26arXi}, with the diagonal $\kappa_k=\frac1n\bigl(p-F(k)+\tau(1-\mu^{pk})\bigr)$ in place of the $\kappa_k$ there. Simplicity of $\kappa_1$ follows as in \cite[Lemma~12]{Che26arXi}. For $j=n-1$ the numerator vanishes.

For $j=1,$ use $1-\mu^{2p}=(1-\mu^p)(1+\mu^p),$ $F(2)-F(1)=\frac{(1-\mu^p)(\mu^p-\mu)}{1-\mu^2}$ and $\mu^{2p}-\mu^p=-\mu^p(1-\mu^p)$. The quotient becomes $(1+\mu^p)/b(\tau),$ and $1+\mu^p=2\cos\frac{\pi p}n\,e^{-i\pi p/n},$ so
\[
c_1(\tau)=\frac2n\sin\frac{2\pi p}{n}\;\Imm\frac{e^{-i\pi p/n}}{b(\tau)}=-\frac2n\sin\frac{2\pi p}{n}\;\frac{\Imm\bigl(b(\tau)e^{i\pi p/n}\bigr)}{|b(\tau)|^2}.
\]
Since $\frac{\mu^p-\mu}{1-\mu^2}=-\frac{\sin(\pi m/n)}{\sin(2\pi/n)}\,e^{-i\pi m/n},$ we get $\Imm\bigl(b(\tau)e^{i\pi p/n}\bigr)=\tau\sin\frac{\pi p}n-\frac{\sin(\pi m/n)}{2\cos(\pi/n)},$ which gives \eqref{e:c1}. The two expressions for $\tau^*_m$ agree because $\sin\frac{\pi m}n+\sin\frac{\pi(m+2)}n=2\cos\frac\pi n\sin\frac{\pi(m+1)}n$.

For even $n,m$ and $j=n/2$ we have $\mu^{j+1}=-\mu$ and $p$ odd. The numerator is $1+\mu^p,$ and the denominator is $\frac{1+\mu^p}{1+\mu}-\frac{1-\mu^p}{1-\mu}-2\tau\mu^p=2b(\tau)$.
\end{proof}

At $\tau=0$ and $\tau=1$, admissible one-sided first-harmonic perturbations are obtained by taking $\eta_k=1+\cos(2\pi k/n)$ and $\eta_k=\cos(2\pi k/n)-1$, respectively.
The added constants affect only $\hat\eta(0)$ and therefore leave the second-order normal coefficient unchanged.

At $\tau=0,$ \eqref{e:c1} reduces to \cite[(14)]{Che26arXi}. Thus the first harmonic pushes the eigenvalue outward exactly for $\tau<\tau^*_m$. For $m=1$ we get $\tau^*_1=1/(4\cos^2\frac\pi n),$ for instance $\tau^*_1=\frac{3-\sqrt5}2$ for $n=5,$ $\frac13$ for $n=6,$ $1-\frac1{\sqrt2}$ for $n=8,$ and $\tau^*_1\to\frac14$ as $n\to\infty$.

\begin{proposition}\label{p:five}
For $n=5,$ $m=1,$ the chord perturbations push $\lambda_1+\tau g_2$ outward at second order if and only if $\tau<\tau^*_1=\frac{3-\sqrt5}2$.
\end{proposition}

\begin{proof}
The relevant coefficients are $c_1$, given by \eqref{e:c1}, $c_4\equiv0$, and $c_2+c_3$, the coefficient of $|\hat\eta(2)|^2=|\hat\eta(3)|^2$. Each $c_j(\tau)$ has the form $\Imm\bigl(A/(B+\tau C)\bigr)$ with $B+\tau C\ne0$ on $[0,1]$, so its sign is the sign of a linear function of $\tau$, and it suffices to evaluate it at $\tau=0$ and $\tau=1$.
A direct computation gives $c_2(0)=c_2(1)=-\frac{\sqrt5}{10}$, $c_3(0)=-\frac{5+\sqrt5}{20}$ and $c_3(1)=-\frac{\sqrt5}{10}$. 
Hence $c_2<0$ and $c_3<0$, and so $c_2+c_3<0$, on $[0,1]$.
\end{proof}

\begin{conjecture}\label{c:harm}
For all $n\ge5$ and $2(m+1)<n,$ every pair sum $c_j(\tau)+c_{n-j}(\tau),$ $2\le j\le n/2,$ is nonpositive for $\tau^*_m\le\tau\le1$.
\end{conjecture}

We checked Conjecture~\ref{c:harm} numerically for all $n\le30$. Individually, the $c_j$ with $j\ge2$ can be positive beyond $\tau^*_m,$ but the pair sums are not. By the case $j=n/2$ of Theorem~\ref{t:thresh}, for even $n$ and $m$ the alternating harmonic switches off at the same $\tau^*_m$.

\begin{remark}[the threshold point]\label{r:Q}
The point $q_m=\lambda_m+\tau^*_mg_{m+1}$ equals $\lambda_m+\frac{\mu(1-\mu^m)}{n(1+\mu)};$ this is the computation of \cite[Proposition~19(ii)]{Che26arXi} with $2l$ replaced by $m$. Hence $\tau^*_{2l}=t_l$ in the notation of \cite[(15)]{Che26arXi}. For even $n$ and even $m=2l,$ the point $q_m$ is the vertex $\lambda^{(n/2)}_l$ of $S_{n/2},$ which lies in the relative interior of $[\lambda_m,\lambda_{m+1}]$.
\end{remark}

\section{Alternating chords for even $n$ and $m$}

Let $n$ and $m$ be even, let $p=m+1$, $w\in[0,1]$, and let $\omega(w)_k=w$ for even $k$ and $0$ for odd $k$. 
Let $2(m+1)<n$: this condition is needed for the arc to leave $S_n$.

\begin{proposition}\label{p:alt}
The subspace $\operatorname{span}\{f_1,f_{1+n/2}\}$ is invariant under $N_m(\omega(w))$. In the basis $(f_1,$ $f_{1+n/2}),$ the restriction of\, $\frac1nN_m(\omega(w))$ is
\[
\begin{pmatrix}\kappa_a+\frac w2\alpha&\frac w2\beta\\[2pt] \frac w2\alpha&\kappa_b+\frac w2\beta\end{pmatrix},
\]
where
\[
\kappa_a=\lambda_m,\ \ \kappa_b=\frac1n\sum_{s=1}^m\bigl(1-(-\mu)^s\bigr),\ \ \alpha=\frac{1-\mu^p}{n},\ \ \beta=\frac{1+\mu^p}{n},
\]
so the eigenvalue $\lambda(w)$ continuing $\lambda_m$ is a root of an explicit quadratic equation. At $w=1,$ the vector $x$ with $x_{2i}=x_{2i+1}=\mu^{2i}$ is an eigenvector with eigenvalue $\lambda^{(n/2)}_{m/2}=\frac2n\sum_{l=1}^{m/2}(1-\mu^{2l})=q_m$.
\end{proposition}

\begin{proof}
Write $D_{\omega(w)}=\frac w2(I+J)$ with $J=\operatorname{diag}((-1)^k)$. Then $Jf_1=f_{1+n/2}$ and $Jf_{1+n/2}=f_1$. Moreover, $f_1$ and $f_{1+n/2}$ are eigenvectors of $Q$ with eigenvalues $\mu$ and $-\mu,$ which gives the matrix.

At $w=1,$ even vertices have arcs of lengths $1,\dots,m+1,$ and odd vertices arcs of lengths $1,\dots,m$. For both $2i$ and $2i+1,$ the arc of length $1$ either joins the two members of the pair $\{2i,2i+1\},$ on which $x$ is constant, or is part of the common target set $\{2i+2,\dots,2i+m+1\}$. These targets form the $m/2$ pairs $\{2i+2l,2i+2l+1\},$ on which $x=\mu^{2(i+l)}$. Hence $(N_mx)_{2i}=(N_mx)_{2i+1}=2\sum_{l=1}^{m/2}(\mu^{2i}-\mu^{2i+2l})$. The last equality is the one from Remark~\ref{r:Q}.
\end{proof}

At $w=1$ the trace of this matrix equals $q_m+\frac{m+1}n$, so its eigenvalues are $q_m$ and the real number $\frac{m+1}n$. In all cases we computed, the two eigenvalues stay at distance at least $0.14$ for $0\le w\le1$, and the branch $\lambda(w)$ starting at $\lambda_m$ ends at $q_m$. 
Thus $w\mapsto\lambda(w)$ is an explicit arc from $\lambda_m$ ($w=0$) to the vertex $q_m$ of $S_{n/2}$ ($w=1$). It leaves $S_n$ at once, with $d_m(\lambda(w))=\frac{c_1(0)}8w^2+O(w^3)$ \cite[Theorem~15]{Che26arXi}. Numerically it stays outside $S_n$ for $0<w<1,$ with maximal distance $9.809\cdot10^{-4}$ ($n=8,$ $m=2,$ at $w\approx0.630$), $6.895\cdot10^{-4}$ ($10,2$), $4.528\cdot10^{-4}$ ($12,2$) and $2.110\cdot10^{-4}$ ($12,4$). At $w=1,$ as in the balanced blow-up \cite[Proposition~18]{Che26arXi}, each pair $\{2i,2i+1\}$ acts as a single vertex, and the quotient is the circulant of order $n/2$ realizing the vertex of $S_{n/2}$.

\begin{conjecture}\label{c:eight}
Over the edge $[\lambda_2,\lambda_3]$ of $S_8,$ the boundary of $\Lambda_8$ is the arc $\{\lambda(w):0\le w\le1\}$ of Proposition~\textup{\ref{p:alt}} followed by the segment $[q_2,\lambda_3];$ the two meet at a corner at the vertex of $S_4$.
\end{conjecture}

\section{The boundary of $\Lambda_5$}\label{s:five}

By the symmetries $z\mapsto\bar z$ and $z\mapsto1-\bar z$ of $\Lambda_5$ \cite[Remark~1]{Che26arXi}, and since $[0,\lambda_1]\subset\partial\Lambda_5$ \cite[Theorems~5, 6]{AgaChe05LAA} (based on \cite{DmitrievDynkin45E,DmitrievDynkin46E}), it suffices to study the nontrivial upper-boundary portion adjacent to $[\lambda_1,\lambda_2]$.

Consider the cycle with chords $k\to k+2$ of weights $(\omega_0,\dots,\omega_4),$ and put $u=1-5\lambda$. Let $P(u;\omega)$ be the characteristic polynomial with the factor for $\lambda=0$ removed. It is multilinear in $\omega$ and has 55 terms. We use the four families
\[
\mathrm A:\ (a,s,a,0,0),\qquad \mathrm{T_1}:\ (a,1,a,0,0),\qquad \mathrm{T_2}:\ (a,1,a,b,b),\qquad \mathrm B:\ (a,s,a,b,b),
\]
which are invariant under the reflection $k\mapsto2-k;$ it fixes $1$ and swaps $0\leftrightarrow2$ and $3\leftrightarrow4$. For A, the characteristic equation \textup(including the root $u=1,$ i.e., $\lambda=0$\textup) reads
\[
(u+a)^2(u+s)u^2=1+2as+(2a+s+a^2s)\,u+a^2u^2,
\]
and for $\mathrm T_1$ the equation factors as $(u+1)\bigl(u^3+(1+2a)u^2+(1+a)^2u+(1+2a)\bigr)=0$.

In each family, $\lambda$ denotes the eigenvalue obtained by continuation from $\lambda_1$ (the root $u=\mu$ at $a=s=b=0$) along the pieces described below; it remains simple along them. All weights lie in $[0,1]$, and each piece is the component, containing the previous transition point, of the set of values of this eigenvalue subject to the indicated equations.
\begin{conjecture}\label{c:five}
The upper boundary of $\Lambda_5$ between $\lambda_1$ and $\lambda_2$ consists of the following pieces\textup, in order from $\lambda_1$\textup:
\begin{enumerate}
\item[(A)] the component through $\lambda_1$ $(a\!=\!s\!=\!0)$ of the curve $\{P=0,\ \Imm(P_a\overline{P_s})=0\}$ for family $\mathrm A,$ up to the point where $s\!=\!1;$ it lies on an algebraic curve of degree $20$ in $5\lambda$\textup, irreducible over $\mathbb Q;$
\item[($\mathrm{T_1}$)] the corresponding eigenvalue curve of family $\mathrm T_1$ \textup(implicit degree $6$\textup{),} up to the point where the switching function of the chords at\/ $3,4$ vanishes\textup;
\item[($\mathrm{T_2}$)] the curve $\{P=0,\ \Imm(P_a\overline{P_b})=0\}$ for $\mathrm T_2,$ up to the point where the switching function of the middle chord vanishes\textup;
\item[(B)] the curve $\{P=0,\ \Imm(P_a\overline{P_b})=\Imm(P_a\overline{P_s})=0\}$ for $\mathrm B,$ which ends tangentially on the edge at $q_1=\lambda_1+\tau^*_1g_2,$ where $a=s=b=\tau^*_1;$
\item[(E)] the segment $[q_1,\lambda_2]$.
\end{enumerate}
\end{conjecture}

\begin{figure}[ht]
\centering
\includegraphics[width=\textwidth]{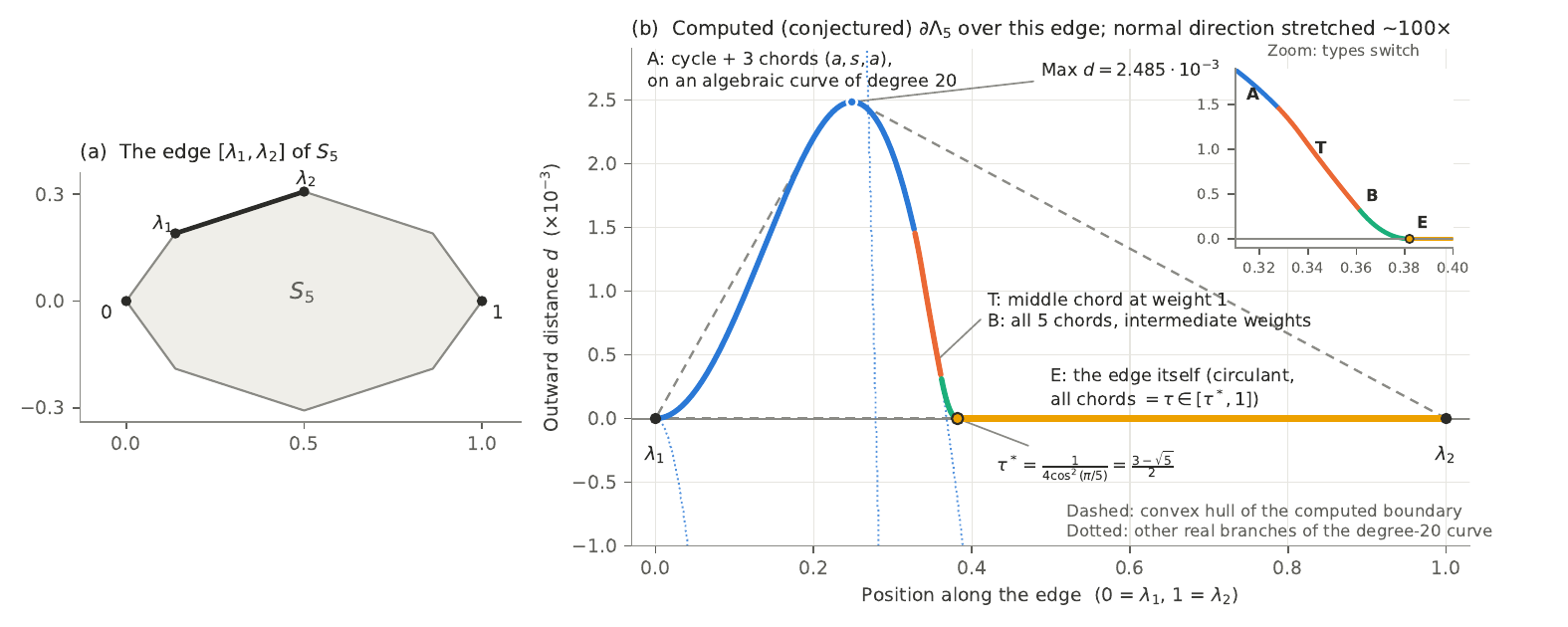}
\caption{Left: the edge $[\lambda_1,\lambda_2]$ of $S_5$. Right: the computed boundary of $\Lambda_5$ over this edge (conjecturally $\partial\Lambda_5$, Conjecture~\ref{c:five}), in edge coordinates with the normal direction stretched about $100$ times. Pieces A, $\mathrm T_1$, $\mathrm T_2$, B and the segment E of Conjecture~\ref{c:five}; dashed: the convex hull of the computed boundary; dotted: the other real branches of the degree-$20$ curve containing piece A.}\label{f:five}
\end{figure}

The transition points were computed at high precision; rounded values are listed in Table~\ref{tab:five}, and the boundary is drawn in Figure~\ref{f:five}.

\begin{table}[ht]
\centering
\small\setlength{\tabcolsep}{3pt}
\caption{Transition points of $\partial\Lambda_5$ over $[\lambda_1,\lambda_2]$. Here $t$ is the parameter along the edge and $d$ the outward distance.}\label{tab:five}
\begin{tabular}{lllll}
\toprule
point & $\lambda$ & $t$ & $d$ & weights\\
\midrule
$\lambda_1$ & $0.138196601125+0.190211303259i$ & $0$ & $0$ & $a=s=0$\\
max of $d$ & $0.227433798077+0.221818990483i$ & $0.24877$ & $2.484887\cdot10^{-3}$ & $a\approx0.3442,$ $s\approx0.7309$\\
$\mathrm A\to\mathrm T_1$ & $0.256132947983+0.230087070158i$ & $0.32723$ & $1.479773\cdot10^{-3}$ & $a=0.478360,$ $s=1$\\
$\mathrm T_1\to\mathrm T_2$ & $0.259886288011+0.230965075262i$ & $0.33733$ & $1.154959\cdot10^{-3}$ & $a=0.506838,$ $b=0$\\
$\mathrm T_2\to\mathrm B$ & $0.268888948044+0.233011324119i$ & $0.36150$ & $3.190825\cdot10^{-4}$ & $a=0.518594,$ $s=1,$ $b=0.036189$\\
$q_1$ & $0.276393202250+0.235114100917i$ & $\frac{3-\sqrt5}2$ & $0$ & $a=s=b=\tau^*_1$\\
\bottomrule
\end{tabular}
\end{table}

\begin{proposition}\label{p:fivecert}
The matrix $\frac15L,$ where $L$ is the Laplacian of the $5$-cycle \textup{(\xz}arc weights $1$\textup) with chords $0\to2,$ $1\to3,$ $2\to4$ of unrounded weights $0.3442,$ $0.7309,$ $0.3442,$ has an eigenvalue $\lambda$ with
\[
d_1(\lambda)=0.0024848865107148434523\ldots,
\]
whose projection onto the line through $[\lambda_1,\lambda_2]$ falls inside the edge\textup, at parameter $0.24876789\ldots$ The analogous four-chord digraph for $n=7$ with weights $(0.1913,0.5884,0.5884,0.1913)$ gives $d_1=0.0012147405919\ldots,$ and the three-chord digraph for $n=6$ with weights $(0.3327,0.739,0.3327)$ gives $d_1=0.0017915891646\ldots$
\end{proposition}

The values were certified in ball arithmetic (Arb via python-flint, $200$ bits) from the exact integer characteristic polynomials, as in \cite[Theorem~5]{Liu026b}. They improve the explicit examples of \cite{Che26arXi,Liu026b} slightly for $n=5,6,$ and noticeably for $n=7,$ where the best three-chord pattern gives $1.194\cdot10^{-3}$ \cite[Table~1]{Liu026b}.

\emph{Evidence for Conjecture~\textup{\ref{c:five}.}}
\begin{itemize}
\item The pieces were found by maximizing the distance from $\frac12$ along rays, over all $20$ off-diagonal entries. For every ray we kept the best of forward and backward continuation sweeps, which were repeated until stable.
\item At the point of maximal distance, the switching law holds with strict complementarity: $|G_{ij}|\ge4.8\cdot10^{-3}$ on all arcs of weight $0$ or $\frac1n,$ while $|G_{ij}|<10^{-8}$ on the arcs with intermediate weights. The transition points of Table~\ref{tab:five} are exactly the points where one switching function vanishes.
\item As an independent check, we maximized $d_1$ over all of $\cL_5$ in $15$ windows of the edge parameter, with $60$ random starts each. Nothing beyond the curve was found, apart from overshoots $\le10^{-5}$ caused by the window penalty.
\item The point $q_1$ is where the second-order analysis predicts the return to the edge (Proposition~\ref{p:five}).
\item The curve of piece A passes through $\lambda_1,$ as it must by the tangential exit \cite[Theorem~9]{Che26arXi}.
\end{itemize}

What is \emph{proved} is the inclusion: 
every point on the selected admissible branches A--B is realized by a standardized Laplacian.
Since $\Lambda_5$ is star-shaped with respect to $\frac12,$ the region between these curves and $[\lambda_1,q_1]$ lies in $\Lambda_5,$ provided the curve is a radial graph about $\frac12,$ as it is numerically. The reverse inclusion, namely that nothing lies beyond the curve, is open (Section~\ref{s:open}).

\section{Other orders}

\begin{figure}[ht]
\centering
\includegraphics[width=\textwidth]{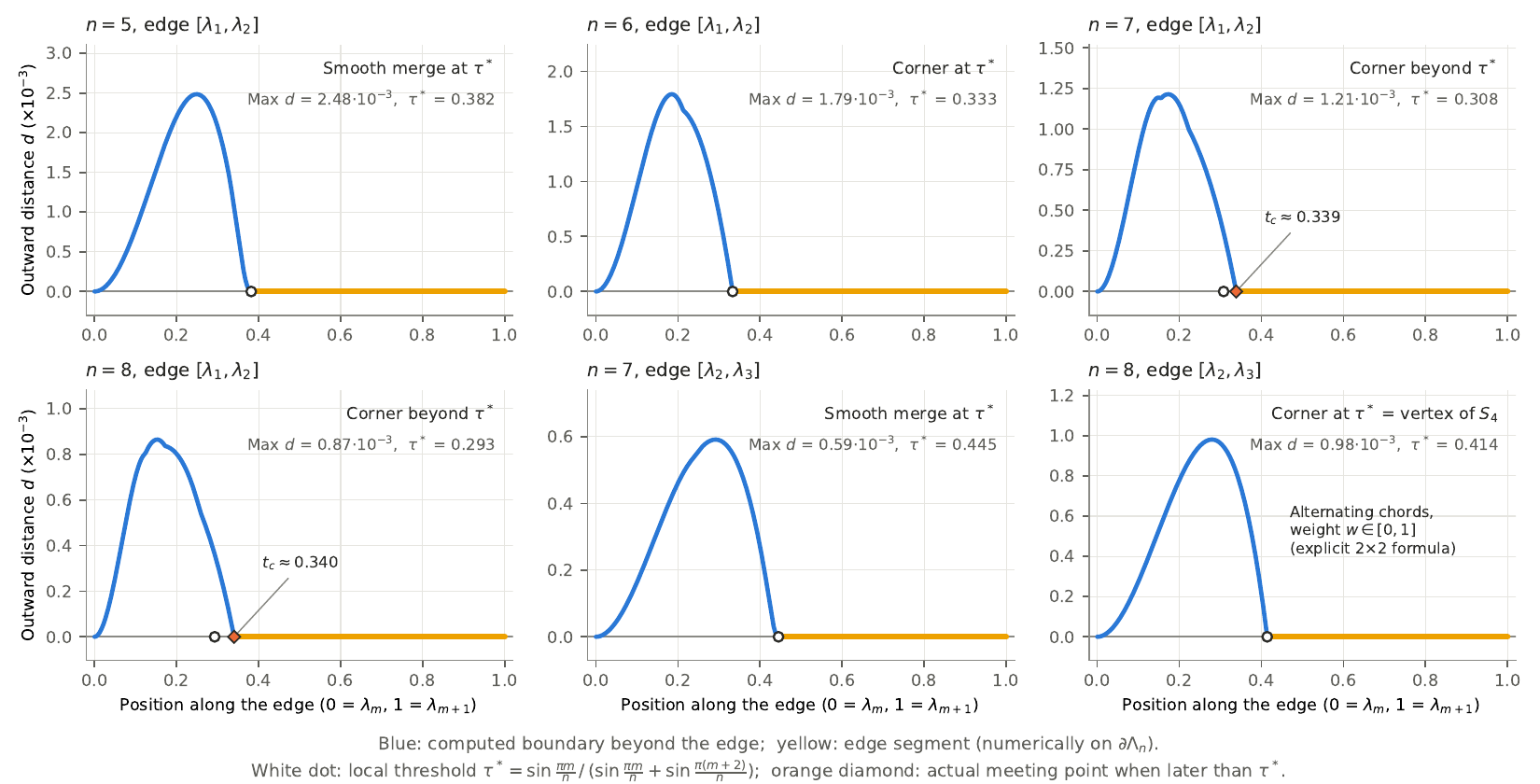}
\caption{The computed boundary of $\Lambda_n$ over the edges $[\lambda_m,\lambda_{m+1}]$ (normal direction stretched). White dots: $\tau^*_m$; orange diamonds: the actual meeting point when it lies beyond $\tau^*_m$.}\label{f:gallery}
\end{figure}

Figure~\ref{f:gallery} shows the boundary over the edges $[\lambda_1,\lambda_2]$ for $n=5,\dots,8$ and over $[\lambda_2,\lambda_3]$ for $n=7,8,$ computed in the same way. 
In all computed cases the cap begins tangentially at $\lambda_m$. For $m=1$ this is a theorem: points of $\Lambda_n$ beyond the edge at distance of order $t^2$ from its line and of order $t$ from $\lambda_1$ exist by \cite[Theorem~15]{Che26arXi}, and $d_1(z)\le C|z-\lambda_1|^2$ for all $z\in\Lambda_n$ near $\lambda_1$ by \cite[Theorem~4]{Liu026b}. For $m\ge2$, only the first of these statements is known \cite[Remark~14]{Liu026b}.
The cap returns to the edge in one of three ways:
\begin{itemize}
\item a \emph{smooth merge at $\tau^*_m$}: $(n,m)=(5,1)$ and $(7,2)$. The chord weights converge to the uniform value $\tau^*_m,$ consistent with Theorem~\ref{t:thresh}.
\item a \emph{corner at $\tau^*_m$}: $(6,1),$ where a four-chord family crosses the edge line at $0.3333\pm10^{-4};$ and $(8,2),$ where the alternating family of Proposition~\ref{p:alt} ends at the vertex of $S_4$.
\item a \emph{corner beyond $\tau^*_m$}: $(7,1)$ and $(8,1),$ with meeting points $t_c\approx0.339$ and $0.340,$ against $\tau^*_1=0.308$ and $0.293$. On $(\tau^*_m,t_c)$ the points outside $S_n$ are realized by matrices far from the circulant. So $\tau^*_m$ is only a lower bound for the return point: the edge cannot be on $\partial\Lambda_n$ before $\tau^*_m$ (by Theorem~\ref{t:thresh} and the argument of \cite[Remark~9]{Che26arXi}), but globally the return may come later.
\end{itemize}

\emph{Growing palindromic blocks.} In the numerically identified boundary families for $m=1$, the realizing digraph consists of the Hamiltonian cycle together with a palindromic block of consecutive chords $k\to k+2$ (cf.\ Lemma~\ref{l:rev}). The block grows outward along the cap:
\begin{itemize}
\item $n=5$: $3\to5$ chords, smoothly, through the saturation of the middle chord;
\item $n=6$: $3\to4$ chords;
\item $n=7$: $3\to4\to5$ chords;
\item $n=8$: $3\to4\to5\to6$ chords.
\end{itemize}
A block of odd length is symmetric about a chord, and a block of even length about the gap between two chords. These are different families, and where they exchange the maximum the boundary has a small corner, visible in Figure~\ref{f:gallery} for $n=6,7,8$. Growth by an outer pair of the same parity is smooth.

\emph{Horizontal edges.} Maximizing $\Imm\lambda$ over all off-diagonal entries (150 random starts) gives $\frac16\cot\frac\pi6$ and $\frac18\cot\frac\pi8,$ both to $5\cdot10^{-16}$. Thus, numerically, the horizontal edges lie on $\partial\Lambda_n$ for $n=6,8,$ in support of an affirmative answer to \cite[Question~3]{Che26arXi} and \cite[Question~19]{Liu026b} in these cases. Together with Figure~\ref{f:gallery}, this gives a conjectural complete picture of $\partial\Lambda_6$ as well (tested less thoroughly than for $n=5$): two caps over $[\lambda_1,\lambda_2]$ and its mirror edge, plus their conjugates, and the edges of $S_6$ elsewhere.

\section{Open problems}\label{s:open}

\begin{enumerate}
\item \emph{An upper bound.} Prove that nothing lies beyond the caps, at least for $n=5,$ or at least the single inequality $\max_{\Lambda_5}d_1=0.0024848865\ldots$
   \begin{itemize}
   \item By the linear-programming argument of \cite[Proposition~17]{Che26arXi}, \cite[Lemma~16]{Liu026b}, it suffices to consider matrices with at most two intermediate weights per row and at most $2n-2=8$ in total.
   \item This is far beyond the six-dimensional certificates of \cite{Liu026b,Marshall026}.
   \item Classical numerical-range bounds \cite{Bendixson02,Pick1922,Browne30} are much too weak: the support function of $\bigcup_{\LT}W(\LT)$ is maximized at $0/1$ matrices and is thus computable, but for $n=5$ in the direction $\nu_1$ it exceeds $h_{S_5}$ by $0.097,$ about $40$ times the height of the cap.
   \end{itemize}
\item \emph{Conjecture}~\ref{c:harm}, and a local theorem: for $\tau^*_m<\tau\le1$, no matrix in $\cL_n$ sufficiently close to $\frac1nN_m(\tau\mathbf 1)$ has an eigenvalue beyond the line of the edge near $\lambda_m+\tau g_{m+1}$. This would be an analogue of \cite[Proposition~3, Lemma~13]{Liu026b} away from the vertex. By the late return observed for $m=1$, $n=7,8$, such a local statement does not imply that the edge point lies on $\partial\Lambda_n$.
\item \emph{The late return for $m=1,$ $n\ge7$}: which family realizes the cap on $(\tau^*_1,t_c),$ and what determines $t_c$?
\item \emph{Convergence to the cycloidal region.} By \cite[Theorem~7]{AgaChe05LAA}, $S_n$ converges in the Hausdorff distance to a cycloidal region $C$. Does $\Lambda_n$ converge to $C$ as well \cite[Question~4]{Che26arXi}? Since $S_n\subseteq\Lambda_n$, this is equivalent to $\sup_{\lambda\in\Lambda_n}\operatorname{dist}(\lambda,S_n)\to0$. Two features of $C$ are already known to be inherited by $\Lambda_n$:
\begin{itemize}
    \item the height: $z(n)=\max\{|\Imm\lambda|:\lambda\in\Lambda_n\}\to\frac1\pi$, since $\frac1n\cot\frac\pi n\le z(n)\le\frac1{2n}\cot\frac\pi{2n}$, by Pick's bound (\cite[Theorem~5]{AgaChe05LAA}, based on \cite{Pick1922}), with equality on the right for odd $n$ \cite[Proposition~2]{AgaChe05LAA};
    \item the tangent directions at $0$ and $1$, by the angles of Dmitriev and Dynkin \cite[Theorem~5]{AgaChe05LAA}.
\end{itemize}
    These outer bounds, however, only confine $\Lambda_n$ to a region converging to $[0,1]\times\bigl[-\frac1\pi,\frac1\pi\bigr]$.
\item \emph{The size of the caps.} Numerically, the caps have height of order $n^{-3}$ ($n^3h\approx0.3$--$0.6$ for $n\le40$), and the best chord block has length close to $n/2$. This is much smaller than $d_H(S_n,C)$, which is numerically $\approx\frac{\pi}{12}\,n^{-2}$. Determine the asymptotics of the caps. Any bound $o(1)$ on $\sup_{\lambda\in\Lambda_n}\operatorname{dist}(\lambda,S_n)$ would answer the previous question affirmatively.
\end{enumerate}

\subsection*{Acknowledgments}
This work was partially supported by the Israel Science Foundation (grant 907/24).
Computations and parts of the analysis were carried out with the help of AI assistants.

\bibliographystyle{abbrvurl}   
\bibliography{all2}

\begin{thebibliography}{10}

\bibitem{AgaChe05LAA}
R.~Agaev and P.~Chebotarev.
\newblock On the spectra of nonsymmetric {Laplacian} matrices.
\newblock {\em Linear Algebra and its Applications}, 399:157--168, 2005.
\newblock \href {https://doi.org/10.1016/j.laa.2004.09.003}
  {\path{doi:10.1016/j.laa.2004.09.003}}.

\bibitem{Bendixson02}
I.~Bendixson.
\newblock Sur les racines d'une {\'e}quation fondamentale.
\newblock {\em Acta mathematica}, 25:359--365, 1902.
\newblock \href {https://doi.org/10.1007/BF02419030}
  {\path{doi:10.1007/BF02419030}}.

\bibitem{Browne30}
E.~T. Browne.
\newblock The characteristic roots of a matrix.
\newblock {\em Bulletin of the American Mathematical Society}, 36:705--710,
  1930.
\newblock \href {https://doi.org/10.1090/S0002-9904-1930-05041-7}
  {\path{doi:10.1090/S0002-9904-1930-05041-7}}.

\bibitem{Che26arXi}
P.~Chebotarev.
\newblock Digraph {Laplacian} spectra: The 2005 polygon conjecture fails for
  all orders above four.
\newblock arXiv:2610.14100v3 [math.CO], 2026.
\newblock \href {https://doi.org/10.48550/arXiv.2610.14100}
  {\path{doi:10.48550/arXiv.2610.14100}}.

\bibitem{DmitrievDynkin45E}
N.~Dmitriev and E.~Dynkin.
\newblock On the characteristic numbers of a stochastic matrix.
\newblock {\em Dokl. Akad. Nauk SSSR}, 49(3):159--162, 1945.
\newblock In Russian.

\bibitem{DmitrievDynkin46E}
N.~Dmitriev and E.~Dynkin.
\newblock Characteristic roots of stochastic matrices.
\newblock {\em Izv. Akad. Nauk SSSR, Ser. Mat.}, 10:167--184, 1946.
\newblock In Russian. Translated in: American Mathematical Society Translations
  (2), 140 (1988), 57--77, doi:10.1090/trans2/140/04.

\bibitem{Karpelevich51E}
F.~I. Karpelevi\v{c}.
\newblock On the characteristic roots of matrices with nonnegative elements.
\newblock {\em Izv. Akad. Nauk SSSR, Ser. Mat.}, 15:361--383, 1951.
\newblock In Russian. Translated in: American Mathematical Society Translations
  (2), 140 (1988), 79--100, doi:10.1090/trans2/140/05.

\bibitem{Liu026}
D.~Z. Liu.
\newblock The {Agaev--Chebotarev} polygon conjecture holds exactly for $n\le4$.
\newblock Preprint, Zenodo, posted 25 October, 2026.
\newblock \href {https://doi.org/10.5281/zenodo.23630989}
  {\path{doi:10.5281/zenodo.23630989}}.

\bibitem{Liu026b}
D.~Z. Liu.
\newblock Non-convexity of the eigenvalue region of standardized {Laplacians}
  and the case $n=4$ of the {Agaev--Chebotarev} conjecture.
\newblock Preprint, Zenodo, posted 30 October, 2026.
\newblock \href {https://doi.org/10.5281/zenodo.23812972}
  {\path{doi:10.5281/zenodo.23812972}}.

\bibitem{Marshall026}
L.~Marshall.
\newblock ai-math-results: {AI}-generated, computer-checked mathematical
  results, folder results/polygon-n4.
\newblock Zenodo, release v2026.10.29c (all versions:
  doi:10.5281/zenodo.23748157), 2026.
\newblock \href {https://doi.org/10.5281/zenodo.23761284}
  {\path{doi:10.5281/zenodo.23761284}}.

\bibitem{Pick1922}
G.~Pick.
\newblock {\"Uber die Wurzeln der charakteristischen Gleichungen von
  Schwingungs\-problemen}.
\newblock {\em Zeitschrift f\"ur angewandte Mathematik und Mechanik},
  2:353--357, 1922.
\newblock \href {https://doi.org/10.1002/zamm.19220020503}
  {\path{doi:10.1002/zamm.19220020503}}.

\end{thebibliography}

\end{document}